\magnification=\magstep1
\def\Z{{\bf Z}}
\def\C{{\cal C}}
\def\Q{{\bf Q/Z}}
\def\N{{\bf N^{+}}}
\def\NN{{\bf N^*}}
\def\z{{\bf Z_{\lambda}}}
\def\l{\lambda}
\def\L{\Lambda}
\def\Z{{\bf Z}}

\centerline{\bf An approach to $F_1$ via the theory of $\L$-rings}
\medskip

\centerline{\bf by Stanislaw Betley}

\vskip2cm

{\bf 0. Introduction.}

\bigskip

This  note  is devoted to the preliminary study of the concept of Borger from [B],
that the decent data from $\Z$ to $F_1$ for commutative rings is the $\L$-structure. More precisely, he claims, that  the category of rings over $F_1$ should consist of $\L$-rings and the restriction of scalars from $\Z$ to $F_1$ takes any commutative ring $R$ to its ring of Witt vectors $W(R)$ with the canonical $\L$-structure. In this approach the mythical field $F_1$ is equal to the ring of integers $\Z$ with the canonical $\L$-structure. We will denote it as $\z$ or just  $F_1$ throughout the sections 1-4.

In [Be] we proved that the categorical $\zeta$-function of the category of commutative monoids calculates the Riemann $\zeta$-function of the integers. This was done in order to support the idea, that while trying to view $\Z$ as a variety over the field with one element we should consider integers as multiplicative monoid without addition. The idea that $\Z$ treated as a variety over $F_1$
should live in the category of monoids is well described in the
literature, see for example [KOW] or [D]. But, because category of monoids is too rigid, most authors instead
of working  with monoids directly extend their field of scalars
from $F_1$ to $\Z$ (or other rings), assuming that scalar
extension from $F_1$ to $\Z$ should take a monoid $A$ to its
monoid ring $\Z [A]$. This agrees  well with the expectation that
rings should be treated in the category of monoids as monoids with
ring multiplication as a monoidal operation. Then the forgetful
functor from rings to monoids and the scalar extension as
described above  give us the nice pair of adjoint functors. But
this approach carries one disadvantage. It takes us quickly from
formally new approach via monoids to the classical world of rings
and modules over them  or to other abelian categories. Working with $\L$-rings seems to be a good way of overcoming this disadvantage. We will preserve monoidal point of view but our approach will allow to use more algebra-geometric methods then the ones which are at hand in the world of monoids.

The definition of $F_1$ as above does describe  this object almost like a field. If $\L$-operations are part of the structure then the ideals in our rings should be preserved by them. It is easy to observe that in $\z$ there are no proper ideals preserved by $\L$-operations. We will try to justify the following point of view that in the category of $\L$-rings we can calculate the Riemann $\zeta$ of the integers in two ways. One - as the $\zeta$-function of the category of modules over $\z$ and the other (geometric) via calculating numbers of fixed points of the action of the Frobenius morphism on the affine line over the algebraic closure of $\z$.

In order to show that our program works we will try first of all  to answer the question what is the algebraic closure of $F_1$. Then we show that the $\zeta$-function of the affine line over $F_1$ calculated by the generating function is the same as the $\zeta$-function of integers. Next we have to find the proper category of modules over $F_1$ and calculate its categorical $\zeta$ function.  At the end we prove that the categorical $\zeta$-function of the category of modules over $F_1$ is equal to the $\zeta$-function of the integers.

\bigskip

{\bf I. $F_1$ and the category of commutative monoids.}

\bigskip

Let  $M_{ab}$ denote the category of commutative monoids with unit and unital maps. This section is devoted only to the rough description of  certain features of $M_{ab}$ which are crucial for our approach in the next sections. They are  either obvious or well described elsewhere so we are very brief here.

Let us start from a some piece of notation. We will denote by $\NN$ the multiplicative monoid of natural numbers. The symbol $\N$ denotes the  monoid of natural numbers with addition. Observe that any monoid $M\in M_{ab}$ carries natural action of $\NN$ by identifying $k\in \NN$ with $\psi^k:M\to M$ where $\psi^k(m)=m^k$. This structure is obvious in $M_{ab}$ and adds very little while studying monoids, but should be reflected always, when we want to induce structures from $M_{ab}$ to other (abelian) categories. This action will be addressed  as an action of Adams operations on $M$. We can also interpret it as an action of the powers of the Frobenius endomorphism. Observe that in characteristic $p$ the Frobenius endomorphism acts by rising an element to the $p$th power. We propose to view the Frobenius action in a uniform, characteristic free  way.  In any structure by Frobenius action we mean an action of $\NN$ where $k\in \NN$ acts by rising an element to the $k$th power if this gives a morphism in the considered structure or $k$ acts as identity.

\bigskip

In [D] Deitmar developed  the algebraic geometry in the setting of monoids, associating to $M\in Obj(M_{ab})$ a topological space with a sheaf of monoids which resembles a spectrum of prime ideals of a commutative ring with its structural sheaf. In his language $F_1$ is the one element  monoid consisting only of the unit.

\bigskip
{\bf Definition 1.1:} For any unital monoid $M$ we define the polynomial ring over $M$ by
$$M[X]=M\times \N.$$ We write $ax^n$ for the element $(a,n)$ with
convention $(a,0)=a$. The evaluation of $ax^n$ at the point $b\in
M$ equals $ab^n$.

\bigskip

{\bf Remark 1.2:} Our definition agrees with the general
expectation ( compare [D], [S]) that scalar extension from $F_1$
to $\Z$ should take any monoid $A$ to its integral monoid ring
$\Z[A]$. Indeed we have an isomorphism of rings
$$\Z [M\times \N]\to \Z[M][X]$$
\noindent given by the formula
 $$\mathop\Sigma_{i\in I}z_i(a_i,n_i)\mapsto \Sigma z_ia_iX^{n_i}$$
\bigskip

Since we know what is the polynomial ring over a monoid we can try to imagine what is the algebraic closure of $F_1$.Typically, when we want to close algebraically a field $K$ we have to add to it, at least,  all roots of elements of $K$ (this is sufficient in the finite field case). In other words we can say that $\bar K$ is the minimal field which contains all roots of the elements of $K$ and has no algebraic extensions. This second condition is equivalent to saying that all roots of elements of $\bar K$ are contained in $\bar K$. So it is not surprising to say:
\bigskip

{\bf Remark 1.3:} $\bar F_1= \Q$

\bigskip

For any monoid $M$ we write $M_+$ for $M$ with $0$ added to it. This is important when we want to talk about geometric points of varieties over $F_1$. In [D] prime ideals in $M$ are given by subsets of $M$ satisfying obvious for ideals conditions. Empty set is a good prime ideal, and it is necessary to consider it  at least in order to have one point in  $Spec(F_1)$. When we have a  variety $M$ over $F_1$ (a monoid) and $B$ is an
extension of $F_1$ we can talk about $M(B)$, the $B$-points of
$M$. They are equal(following [D]) to $Hom_{M_{ab}}(M_+,B_+)$. Here we have to add zeros to monoids because a map cannot have an empty value on some element.

We can calculate  the number of points of $\N$, treated as an affine line over $F_1$, in $\bar F_1$. The Frobenius action on $\Q$ extends to the action on $\Q$-points of $\N$. We easily calculate that for $k\in \NN$ we have the following formula for the order of the set of fixed points:
$$\vert (\N(\Q))^k \vert =k$$

If we view the Weil zeta function for
varieties over $F_p$ as a way of keeping track of the numbers of
their points over $\bar F_p$ fixed by the action of the powers of the Frobenius morphisms then we get the  formula for the Riemann zeta function $\zeta_R$ in the similar spirit in $M_{ab}$:

$$\zeta(\N,s)=\mathop\Sigma_{k=1}^{\infty}\vert (\N(\Q))^k \vert^{-s}$$

\bigskip

As we said in the introduction the geometry in $M_{ab}$ is too weak to expect that we can fulfill Delign\'e's program there  and approach the Riemann conjecture. But the picture described above is very basic and leads to the correct $\zeta$ function. This implies that we want to see our more sophisticated structures of $F_1$-objects as the structures induced from the picture described above.

\bigskip

{\bf II.  Preliminaries on $\L$-rings.}
\bigskip
Our rings are always commutative with units. Following [Y, Definition 1.10] we have:
\medskip
{\bf Definition 2.1:} A $\l$-ring is a ring $R$ together with functions $$\l^n:R\to R\ \ \ (n\geq 0)$$ satisfying for any $x,y\in R$:

(1) $\l^0(x)=1$,

(2) $\l^1(x)=x$,

(3) $\l^n(1)=0$ for $n\geq 2$,

(4)$\l^n(x+y)=\sum_{i+j=n}\l^i(x)\l^j(y)$,

(5) $\l^n(xy)=P_n(\l^1(x),...,\l^n(x);\l^1(y),...,\l^n(y))$,

(6) $\l^n(\l^m(x))=P_{n,m}(\l^1(x),...,\l^{nm}(x))$.

\bigskip
Above $P_n$ and $P_{n,m}$ are certain universal polynomials with integer coefficients obtained via symmetric functions theory (see [Y, Example 1.9 and 1.7]). By a homomorphism of $\L$-rings we mean a ring homomorphism which commutes with $\L$-operations. We say that $x\in R$ is of degree $k$, if $k$ is the largest integer for which $\l^k(x)\neq 0$. If such finite $k$ does not exist we say that $x$ is of infinite degree. Observe that (by formula $4$) the map
$$R\ni x\mapsto \mathop\Sigma_{i\geq 0} \l^{i} (x) t^i$$ is a homomorphism from the additive group of $R$ to the multiplicative group of power series over $R$ with constant term $1$. We will denote this map as $\l_t(x)$. Observe also that $\l_t(0)=1$ and hence $\l_t(-r)=\l_t(r)^{-1}$.

Ring of integers $\Z$ carries the unique, canonical $\L$-ring structure  described by the formula
$\lambda^n(m)={m\choose n}$. Similarly all integral monoid rings $\Z[M]$ will be considered with the
$\L$ structure defined for any $m\in M$ by formulas $$\l^1(m)=m$$  $$\l^i(m)=0\ \ {\rm for} \ \ i>1.$$ We will always considered integral monoid rings with such $\L$-structure, because it corresponds well with the monoidal point of view on the field with one elements.
\bigskip
{\bf Lemma 2.2:} Let $R$ be equal to the  monoidal $\L$-ring  ${\Z}[M]$ with the $\L$-structure defined above. Then in $R$  only generators $m\in M\subset {\Z}[M]$ are of degree $1$.
\medskip
Proof.  We know that the elements of $M$ are of degree 1. By separating positive and negative coefficients we get  ${\Z}[M]\ni r=\mathop\Sigma_{i=1}^k a_im_i -\mathop\Sigma_{j=1}^l b_jm_j$, where all $a_i$s and $b_j$s belong to $\Z$ and are greater than $0$.  Observe that the assumption that $r$ is of degree $1$ implies that $\l_t(r)=1+bt$ for a certain $b\in {\Z}[M]$. We can easily calculate $\l_t$-functions in the case of monoidal rings. Let $m\in M$ and a be a positive integer. Then
$$\l_t(am)=\l_t(m)^a=(1+mt)^a$$
Hence we easily get

$$ \l_t(\mathop\Sigma_{i=1}^k a_im_i-\mathop\Sigma_{j=1}^l b_jm_j)= \prod_{i=1}^k(1+m_it)^{a_i}/\prod_{j=1}^l(1+m_jt)^{b_j}$$
If $r$ is of degree $1$ we have equality

$${\bf (*)}\ \ \ \prod_{i=1}^k(1+m_it)^{a_i}=(1+bt)\prod_{j=1}^l(1+m_jt)^{b_j}$$

From this, by comparing coefficients at the highest degree of $t$ we get
$$b=\prod_{i=1}^k(m_i)^{a_i}/\prod_{j=1}^l(m_j)^{b_j}$$
and hence $b\in M$. On the other hand, when we calculate the coefficient at the first degree in the equality {\bf (*)} we get
$$  \mathop\Sigma_{i=1}^k a_im_i=b+ \mathop\Sigma_{j=1}^l b_jm_j$$

But by the definition of $a_i$s and $b_j$s this is possible only when $r=b\in M$.

\bigskip
Let $M_{ab}$ denote the category of commutative monoids with unit and $Ring^{\l}$ stands for the category of commutative unital rings  with $\L$-structure. We have:
\bigskip

{\bf Proposition 2.3:} The functor $M_{ab}\to Ring^{\l}$, which takes a monoid $M$ to  the  $\L$-ring ${\Z}[M]$ with the $\L$-structure defined above has a right adjoint $Ring^{\l}\to M_{ab}$ which takes a $\L$-ring $R$ to the multiplicative monoid $R_1$ of its elements of degree not exceeding $1$.

\medskip
Proof. By  [Y, Proposition 1.13] we know that  in any $\L$-ring the product of 1-dimensional elements is again 1-dimensional (or equal to 0). Hence $R_1$ is a well defined multiplicative submonoid of $R$ considered here as the  multiplicative monoid. If $f\in Mor_{Ring^{\l}}(R,S)$ then $f$ carries $1$ dimensional elements of $R$ to $1$ dimensional elements of $S$ or to $0$ by the definition of a $\L$-homomorphism. Hence our right adjoint is well defined. The rest of the proof is obvious.

\bigskip
The $\L$-operations on a ring $R$ define on it the sequence of Adams operations $\psi^k:R\to R$ which are natural ring homomorphisms. They can be defined by the Newton formula:
$$\psi^k(x)-\l^1(x)\psi^{k-1}(x)+...+(-1)^{k-1}\l^{k-1}(x)\psi^1(x)=(-1)^{k-1}k\l^k(x)$$
For their properties see [Y, chapter 3]. It is straightforward to check that the canonical $\l$-structure on $\Z$   defines trivial Adams operations and the formula  $m\mapsto m^k$ for $m\in M$ determines  the $k$th Adams operation on the monoidal ring $\Z[M]$. Adams operations can be viewed always as an action on a considered structure by the  multiplicative monoid $\NN$ of natural numbers. Every object $M$ of $M_{ab}$ has naturally such a structure as was described in Section 1.  So proposition 2.3 can be viewed as a statement about adjoint functors between categories with objects carrying the action of $\NN$. It is easy to check that the $\NN$ action on  ${\Z}[M]$ given by $k(m)=m^k$ while treated as the action of Adams operations forces to have  $\l$-structure on  ${\Z}[M]$ satysfying  $\lambda^i(m)=0$ for $i>1$.

Observe that if we consider natural numbers $\N$ as a monoid with addition, then for any ring $R$ we have a description of the polynomial ring over $R$  via the formula
$$R[x]=R[\N]=R\otimes {\Z}[\N].$$
Hence we will consider  $Z[\N]$ as  polynomial ring over $F_1$ in the rest of the paper, with $\L$-structure defined like for any other monoidal ring. Moreover, for any $\L$-ring $R$ we have well defined   $\L$-structure on $R[x]$ because tensor product of rings inherits it from the $\L$-structures of the factors.

\bigskip
{\bf Definition 2.4:} Let $R$ be a $\L$-ring and $I$ is an ideal in $R$. We will call it a $\L$-ideal if it is preserved under the action of $\l^k$, for any $k>0$.

\bigskip
It is straightforward to check that  if we divide a $\L$-ring by a $\L$-ideal then $R/I$ carries the induced $\L$-structure and the quotient homomorphism $R\to R/I$ is a homomorphism of $\L$-rings. Of course the opposite is also true: a kernel of the $\L$-rings homomorphism is a $\L$-ideal. For computations  important is that an ideal $I$  in a $\L$-ring $R$ with $\Z$-torsion free quotient is a $\L$-ideal if and only if it is preserved by the Adams operations (see [Y, Corollary 3.16]).

\bigskip

{\bf III.  Algebraic extensions of $F_1$ and its algebraic closure $\bar F_1$.}

\bigskip

As we said in the introduction, the $\L$-ring $\z$ (our hypothetical $F_1$) can be treated as a field because it contains no proper $\L$-ideals. We have defined the ring of polynomials over $F_1$. On the other hand we can always view an algebraic closure of a field via the ring of polynomials and its quotients. If $K$ is a field we know that every algebraic extension of $K$ should be contained in $\bar K$. We know that  every algebraic extension of $K$ is build out of simple extensions $K\subset K(a)$ where $a$ is a root of a non-decomposable polynomial $f_a \in K[x]$. Moreover, for simple extensions we have the formula $K(a)=K[x]/(f_a)$. It means that we can view $\bar K$ as a sum of the fields of the form $L[x]/(f)$ where $L$ is an algebraic extension of $K$ or even as a sum of simple extensions $K[x]/(f)$. We will try to use this point of view in our context but remembering all the time about the $\L$-structures on our objects. Let us start from the following lemma:
\bigskip

{\bf Lemma 3.1:} Every monic  polynomial $f$, which generates the principal $\L$-ideal in $F_1[x]$ has only roots of unity or $0$ as his roots in $\bf C$. If $f(\mu_n)=0$, where $\mu_n$ is  the prime root of $1$ of degree $n$, then $x^n-1\ \vert \ f$ (and hence $n\leq deg(f)$).

\medskip
Proof. Because our polynomials  are monic the quotients of $F_1[x]$ by ideals generated by them are torsion free. Hence instead of working with $\l$-operations we can assume that our ideals are preserved by all  Adams operations. Assume that an ideal $I$ is generated by a polynomial $f(x)=a_nx^n+a_{n-1}x^{n-1}+...+a_0$. Then $\psi^k(a_nx^n+a_{n-1}x^{n-1}+...+a_0)=a_nx^{kn}+a_{n-1}x^{k(n-1)}+...+a_0$. This formula implies that if $a$ is a root of $f$ then $a^k$ is also a root of $f$ for any natural $k$, when we calculate roots in the field of complex numbers. To see this observe that the evaluation at $a$ map $e_a:F_1[x]\to {\bf C}$ is a homomorphism of rings. Let $I=ker(e_a)$ then $(f)\subset I$ and because $\psi^k((f))\subset (f)$ then $\psi^k(f)\in (f)\subset I$.
But for a polynomial $f$ the statement that if $a$ is a root of $f$ then $a^k$ is also a root for any natural $k$ implies directly that $a$ is a root of $1$ or $a=0$. Hence all the zeros  of $f$ are  the roots of unity or are equal to $0$. On the other hand if $a$ is a primitive root of $1$ of degree $s$ and $f(a)=0$ then all the other roots of unity of degree $s$ are among the zeros of $f$, again  by the argument with Adams operations. This implies immediately our statement.
\bigskip

{\bf Definition 3.2:} Let $K$ be a $\L$-ring. We will say that $f\in K[x]$ is non-decomposable if the principal ideal $(f)$ is preserved by the $\L$-operations and there is no decomposition $f=f_1\cdot f_2$ such that $deg(f_1)>0$, $deg(f_2)>0$ and both ideals $(f_1)$ and $(f_2)$ are preserved by the $\L$-operations.

\bigskip
 Lemma 3.1 describes what are the   non-decomposable polynomials in $F_1[x]$. In our approach  we would like to take the following definition for an algebraic extension  in the category of $\L$-rings (field extensions of $F_1$).
\bigskip

{\bf Definition 3.3:} Let $F$ be a $\L$-ring without nilpotent elements. We say that:

{\bf i.} $K\supset F$ is a simple  algebraic extension of $F$  if $K$ is isomorphic to  $F[x]/(f)$ for some  non-decomposable $f\in F[x]$ and $K$ has no nontrivial $\L$-quotients except $F$ itself.

{\bf ii.} $K$ is an algebraic extension of $F$ of finite degree if there is a sequence of simple extensions $F\subset K_1 \subset...\subset K_s \subset K$.

{\bf iii.} $K$ is an algebraic extension of $F$ if every $k\in K$ is contained in a finite degree algebraic extension $L_k$ of $F$, which is contained in $K$.
\bigskip

We say that $L$ is an algebraic closure of $K$, $L=\bar K$, if $L$ is an algebraic extension of $K$, for any finite degree algebraic extension $K\hookrightarrow L_1$ we have a factorization $K\hookrightarrow L_1 \hookrightarrow L$ and $L$ is minimal with respect to this property. We care only about monic  $f$'s because otherwise $F_1[x]/(f)$ contains torsion elements which are always nilpotent in any $\L$-structure. Some explanation why we have to accept that a field has a nontrivial quotient is contained  in  Remark 3.5.

\bigskip

{\bf Proposition 3.4:} $\bar F_1= F_1[\Q].$
\medskip
Proof. Observe first that the simple algebraic extension of $F_1$ is isomorphic to $F_1[x]/(x^p-1)$ for some prime number $p$. This follows directly from 3.1 because if our extension is of positive degree and has presentation  as $F_1[x]/(f)$ then 3.1 implies that $x^k-1$ divides $f$ for a certain positive $k$. But then either our extension has the forbidden quotient $F_1[x]/(x^k-1)$ (we will say in the future that $f$ has a forbidden factor) or  $f=x^p-1$ for a prime $p$. So we see  that every simple algebraic extension of $F_1$ is isomorphic to $F_1[\mu_p]$ for a prime number $p$ which can also be interpreted as a group ring $\Z[C_p]$ for a cyclic group $C_p$ of order $p$.

Now assume that $K$ is a simple algebraic extension of $F_1[\mu_n]$. We know that $K$ has a presentation as $F_1[\mu_n][x]/(f)$ for a certain non-decomposable $f\in F_1[\mu_n][x]$. Observe, that  for any  $a\in F_1[\mu_n]$ , $\psi^n(a)\in \Z$. Moreover for any  $g\in F_1[\mu_n][x]$ , $\psi^n(g)\in F_1[x]$ by the description of Adams operations. We are using here first of all the fact that the $F_1[\mu_n][x]=F_1[\mu_n]\mathop\otimes_{\Z}F_1[x]$ and the $\L$-operations on the tensor product are induced from those on the factors. Secondly, Adams operations are fully described by $\L$-operations and they are natural ring homomorphisms.

 So we know that $\psi^n(f)\in F_1[x]$.  We can apply to it the same procedure as in the proof of 3.1 and get that if $a$ is a root of $\psi^n(f)$ then all its powers have also this property. So any root of $\psi^n(f)$ is a root of $1$ or $0$ and the statement of 3.1 is fulfilled for $\psi^n(f)$. Because $f\vert \psi^n(f)$ we know that all roots of $f$ are roots of $1$ of degree not exceeding $nk$, where $k=deg(f)$ (zero can be excluded from our considerations). The polynomial $f$ cannot have any root of $1$ of degree $n$. If this would be the case then, because $x^n-1$ is fully decomposable in $F_1[\mu_n][x]$, $f$ would have a forbidden quotient. Consider first the special case and assume  that  $f$ has a root of prime degree $p$, and $(n,p)=1$. This implies that $x^p-1\vert \psi^n(f)$. But $f\vert \psi^n(f)$ and $f$ and $x^p-1$ have a common root. Moreover $x^p-1$ has no root besides $1$ in $F_1[\mu_n]$. Hence $x^p-1 \vert f$ and this implies that $f=x^p-1$ or $f$ has forbidden quotient.
 Now we can consider the general case. Let $a$ be a root of $f$. We know that $a^{nk}=1$ and $a^n\neq 1$. This implies that $a$ is a root of $1$ of degree $ns$ and $s\vert k$. So $a^{ns}=(a^{s})^n=1$ and hence $a^s$ is a root of $1$ of degree $n$. By the same considerations as in the special case we get now that  $f=x^p-\mu_n^i\in F_1[\mu_n][x]$ and $p=s$ is a prime.

Observe that with $n$ and $p$ as in the special case  we have $$F_1[\mu_n][x]/(x^p-1)=\Z[C_n][C_p]=\Z[C_n\times C_p]=\Z[C_{np}]=F_1[\mu_{np}]$$ In the general case observe that  $a\mu_n$ is of  degree $pn$ and freely generates $F_1[\mu_n][x]/(x^p-\mu_n^i)$ over $\Z$ so we can write:
 $$F_1[\mu_n][x]/(x^p-\mu_n^i)=\Z[C_{np}]=F_1[\mu_{np}]$$
This implies directly our proposition.

\bigskip
{\bf Remark 3.5:} Unlike $F_1$ all field extensions of $F_1$ have nontrivial quotients so they contain ''ideals''. Obviously  $(x-1)\vert (x^n-1)$ and this implies that there is a $\l$-ring map $F_1[\mu_n]\to F_1$ given by augmentation. This is caused by the fact that we induce our structures from the category $M_{ab}$ where one point monoid $\bf 1$ is a zero object. It means that beside the map ${\bf 1}\to M$ for any $M\in M_{ab}$ we have also $M\to {\bf 1}$ and this should be reflected in the category of $\l$-rings. Moreover in $M_{ab}$ we have a notion of a quotient object, the quotient map there commutes always with the action of $\NN$ so we have to allow the existence of quotients of our "fields".

\bigskip
Since we have the algebraic closure of $F_1$ we can look for the ''powers of the Frobenius morphism''. Recall the well known fact from algebraic geometry over finite fields. If $\bar F_p$ is an algebraic closure of $F_P$ and $X$ is a variety over $F_p$ then $$X(F_{p^n})=X(\bar F_p)^{\bf n}$$ Here $X(K)$ denotes the $K$-points of $X$ and $Y^{\bf n}$ denotes the fixed points of the $n$th power of the Frobenius morphism action on $Y$. We would like to have the similar formula over $F_1$. For the points of the variety over different fields we have  universal solution. If $R$ is a $\L$-ring and $X$ is a variety over $F_1$  then
$$X(R) = Mor(R,X)$$ where the morphism set is taken in the category of varieties. In the  affine case, but what we mean $X$ is described by another $\L$-ring $S$ ($X=spec(S)$) we have as usual $$Mor(R,X)=Hom_{\L -rings}(S,R)$$
Affine case $R={\Z}[M]$ with $M$ - monoid is crucial for us and in this case  we can give precise meaning to the superscript $\bf n$. A $\L$-ring structure on $R$ implies that there  is a homomorphism $i:R\to W(R)$,  where $W(R)$ is a ring of big Witt vectors over $R$. For any $n$ we have a Frobenius morphism $f_n:W(R)\to W(R)$. If we wiew $W(R)$ as invertible formal power series in indeterminant $t$ with addition given by series multiplication and multiplication coming from ghost coordinates then $i(r)=1+rt$ for any $r\in R$ of degree one. For such elements
$$f_n(1+rt)=(1+r^nt)=\psi^n(1+rt)$$
Hence we can assume that Frobenius action of $n\in N$ on ${\Z}[M]$ is realized by rising monoidal generators to the $n$th power. In other words it means this action is realized by the Adams operations.

Let now $X={\Z}[x]$. Then $X(\bar F_1)=Hom_{\L -rings}({\Z}[x],F_1[\Q])$. Observe that every $\L$-homomorphism from ${\Z}[x]$ is determined by the image of $x$ which should be contained in elements of degree not exceeding $1$. In the case of monoidal rings $$ Hom_{\L -rings}({\Z}[x],{\Z}[M])= M_+$$ because as an image of $x$ can be taken any element of $M$ or $0$ and by  2.2 that is all. As usual $M_+=M\cup \{0\}$. The Frobenius action is realized by Adams operations. Hence we calculate that the set $X(\bar F_1)^{\bf n}$
consists of the roots of unity of degree ${n-1}$ plus additional element  $0$ so has cardinality $n$. From this we get  the following formula for the ordinary Riemann  $\zeta $-function $\zeta_{R}$:
$$\zeta_{R}=\mathop\Sigma_{n=1}^{\infty} 1/X({\bar F_1})^{\bf n}$$
where $X$ is an affine line over $F_1$ equal to ${\Z}[x]$.

\bigskip

{\bf IV.  Categorical $\zeta$ function over $F_1$.}
\bigskip
Let us start from recalling after Kurokawa (compare [K]) the
definition of the zeta function of a category with $0$. If  $\C$
is a category with $0$ we say that $X\in Ob(\C)$ is simple if for
any object $Y$ the set $Hom_{\C}(X,Y)$ consists only of
monomorphisms and $0$. Let $N(X)$, the norm of $X$, denote the
cardinality of the set $End_{\C}(X)$. We say that an object $X$ is
finite if $N(X)$ is finite. We denote by $P(\C)$ the isomorphism
classes of all finite simple objects of $\C$. Then we define the
zeta function of $\C$ as

$$ \zeta(s,\C)=\prod_{P\in P({\C})}(1-N(P)^{-s})^{-1}$$

In [K] Kurokawa studied the properties of such zeta functions but
for us the crucial is:
\medskip

{\bf Remark 4.1:} Let  $Ab$ denote the category of abelian groups
 and $\zeta_R$ stands for the Riemann zeta function of the
integers. Then
$$\zeta(s, Ab)=\zeta_R$$
\bigskip
The following observation was the starting point for our considerations and it  underlines the role of the category of abelian monoids. Recall that  $M_{ab}$ denote the category of abelian monoids with unit and
unital maps. In [Be] we proved:

\medskip
{\bf Theorem 4.2:} $$\zeta_R \cdot (1-2^{-s})^{-1}=\zeta(s,M_{ab})$$

\bigskip

Hence the category of monoids carries all the information needed for calculating the Riemann ${\zeta}$-function of the integers. We want to look at its categorical calculation  from 4.1 in a slightly different way, which is suitable for generalizations. First of all we underline that we are working in the category of $\Z$-modules. There are good analogs of the category of modules over an object $X$  of an abstract category $\cal C$ which has $0$ and all finite limits. Beck in [Bec] defined them as abelian group objects in the category of objects over $X$ (see also [H, chapter 2] ). As is shown in [H] the category of abelian group objects in the category  of rings over a given ring $R$ is equivalent to the category of $R$-modules, where an $R$-module $X$ defines the square zero extension of $R$ with $X$ as a square-zero ideal. In the case of $R=\Z$ we get, as expected, the  category of abelian groups. An abelian group $X$ corresponds to the square zero extension  ${\Z}\triangleright X$. The finite simple objects in the category of rings over $\Z$ are easily seen to come from the simple abelian groups (finite cyclic groups $C_p$ of prime order $p$). For a given $p$ we  see that $N(C_p)$ is equal to the cardinality of the set $Hom_{Rings/{\Z}}({\Z}[x],{\Z}\triangleright C_p)$. The polynomial ring ${\Z}[x]$ is  treated as a ring over $\Z$ via the map which takes $x$ to $0$. All this means that we have the geometrical method for calculating  the categorical ${\zeta}$ of integers. We just have to  count the ${\Z}\triangleright C_p$-points of the affine line over $\Z$.

 We can perform the calculation as above for any commutative  ring $R$ because finite simple objects in the category of $R$-modules correspond to the maximal ideals $I\subset R$ with finite quotient and one checks immediately that the cardinality of $Hom_{Rings/{R}}({R}[x],{R}\triangleright R/{I})$ is the same as the cardinality of $Hom_{R-mod}(R/I,R/I)$. Moreover this cardinality is the same as the number of elements of the residue field at the closed point corresponding to $I$ in $Spec(R)$ so we are really calculating the classical $\zeta$-function of an affine variety $Spec(R)$.

Observe that we can perform the same calculations in the category $M_{ab}$, where the role of integers is played by the field of one element in the sense of [D]. But this gives us no new insight because in the world of monoids the field of one element is represented by one point monoid $\bf 1$ consisting of $1$ only, so we have equality of categories  $M_{ab}/{\bf 1}=M_{ab}$. The affine line over $\bf 1$ is equal to the monoid of natural numbers with addition.

We want to promote here the point of view that  calculations of categorical  $\zeta$ function should   have geometric meaning. By this we mean that if in a category $\cal C$  we have ''ground field $A$'' related to it such that all objects in $\cal C$ have the $A$-structure  and moreover we know what is the  affine line $A[x]$ in $\cal C$ then categorical $\zeta$-function of $\cal C$ should be defined in a geometrical way. But of course in such generality we can expect to have artificial objects in $\cal C$.  To exclude them we will call a finite simple object   $P$ in $\cal C$  geometrically finite if $n(P)=\vert Mor_{\cal C}(A[x],P)\vert$ is finite.
Then we define the geometrical $\zeta$-function of $\cal C$ by the formula

$$ \zeta_g(s,\C)=\prod_{P\in P'({\C})}(1-n(P)^{-s})^{-1}$$

\noindent where $P'(C)$ denote the set of isomorphism classes of geometrically finite objects of $\cal C$ with $n(P)>1$. Geometrically finite simple object $P$ with $n(P)>1$ will be called as  non-degenerate  objects of $\cal C$. These are the only objects which are meaningful for calculating the geometrical Riemann $\zeta$-function of $\cal C$. Observe that in $M_{ab}$ or $R-mod$ finite simple objects are geometrically finite and hence in this cases $\zeta_g=\zeta$. But this is not the case of the category of modules over a $\l$-ring, as we will see below.

Our aim is to  apply the described above strategy for calculating $\zeta$ and $\zeta_g$ functions for the category of $\L$-modules   over $F_1$. We know what is the polynomial ring over a $\L$-ring, so we can use affine line $F_1[x]$ in our constructions. We show below that the categorical $\zeta$-function of $F_1$ agrees with the  results of chapter 3, where we calculated it via the algebraic closure and the Frobenius action.

First we should describe the category of modules over $F_1$. This is done in full details in [H, chapter 2] for any $\L$-ring. The constructions uses the functor $W$ from unital commutative  rings to $\L$-rings which takes any ring $R$ to its ring of Witt vectors $W(R)$. Originally the functor $W$ was defined for rings with  multiplicative unit.   But the universal polynomials which define addition, multiplication and opposite in $W(R)$ do not use multiplicative unit so using the same formulas one can define the value of $W$ on the non-unital rings.

If $R$ is a $\L$-ring then it comes with the $\L$-ring map $\l_R:R\to W(R)$ which is defined by lambda operations on $R$. More precisely, if $\L(R)$ denote the ring of invertible formal power series over $R$ then $\l_R=E\circ \l_t$ where $E$ is the Artin-Hasse exponential isomorphism of $\L(R)$ and $W(R)$ (see [Y, chapter 4]) and
$$\l_t(r)=\mathop\Sigma_{i=0}^{\infty}\l^i(r)t^i$$

As it is proved in [H] the category of modules over a $\L$-ring $R$, by which we mean the category $(Ring^{\l}/R)^+$ of abelian group objects in $Ring^{\l}/R$, is equivalent to the category $R-mod^{\l}$ of $\L$-modules over $R$.  A $\L$-module over $R$ is an $R$ module $M$ with a map $\l_M:M\to W(M)$ which is equivariant with respect to the $\L$-structure of $R$. Here $W(X)$ denotes the Witt ring construction applied to the non-unital ring $X$ with trivial multiplication.  It is easy to check that in this case  $W(M)$ has also trivial multiplication  and additively is equal to the infinite product of $M$. It is shown in [H, Lemma 2.2] that we have an isomorphism of rings
$$i: W(R)\triangleright W(M)\to W(R\triangleright M)$$
\noindent which is induced by the canonical inclusions of $R$ and $M$ into $R\triangleright M$.
A $\L$-module $M$ corresponds in the equivalence of $(Ring^{\l}/R)^+$ and $R-mod^{\l}$ to the $\L$-ring $R\triangleright M$ with the $\L$-ring structure defined by the composition

$$R\triangleright M\buildrel{\l_R \oplus \l_M}\over \longrightarrow W(R)\triangleright W(M)\buildrel{i} \over \longrightarrow W(R\triangleright M)$$

We have another description of the category $R-mod^{\l}$ (see [H, Remark 2.6]). If $M$ is an object of this category and  $\l_M:M\to W(M)$ is a structural map then it has components  $\l_{M,n}:M\to M$ because as sets $W(M)=\prod_{N}M$. Easy calculation shows  $\l_{M,n}$ is $\psi_{R,n}$ equivariant, where $\psi_{R,n}$ is the $n$th Adams operation of $R$. This gives us description of the category $R-mod^{\l}$ as a category of left modules over a twisted monoid algebra $R^{\psi}[N]$ where the multiplicative monoid $N$ acts on any object $M$ through the maps $\l_{M,n}$.

With the understanding of the category $R-mod^{\l}$ presented above we can come back to our situation and analyze   the category $F_1-mod^{\l}$. Observe that the Newton formula which relates Adams and $\l$-operations implies  $\psi_{F_1,n}=id$ for any natural $n$. This  implies  that  $\l_{M,1}=id$ and for $n>1$,  $\l_{M,n}:M\to M$ is any (additive) group homomorphism. So we have:

\bigskip
{\bf Lemma 4.3:} Every object  $(M,\l_M)$ in $F_1-mod^{\l}$ consists of an abelian group $M$ and a sequence of group homomorphisms $\l_{M,n}:M\to M$ satisfying  $\l_{M,n}\circ \l_{M,m}= \l_{M,mn}$ and $\l_{M,1}=id$. Morphisms $(M,\l_M)\to (P,\l_P)$  are given as group homomorphisms $f:M\to P$ which satisfy $f\circ \l_{M,n}=\l_{P,n}\circ f$ for any natural $n$. An object  $(M,\l_M)$ is  simple and finite if $M$ is a cyclic group of prime order.
\medskip
Proof. The description of $F_1-mod^{\l}$ was achieved before the statement of the lemma. Observe that an endomorphism $f:M\to M$ of a cyclic group satisfies: for any subgroup $P<M$, $f(P)\subset P$. This implies immediately, that simple objects are as described.

\bigskip
{\bf Lemma 4.4:} Every non-degenerate object in  $F_1-mod^{\l}$ is of the form $(C_p,\l_{C_p})$, where $\l_{C_p,n}=0$ for $n>1$.
\medskip

Proof. First of all, every non-degenerate object $(M,\l_{M})$ in $F_1-mod^{\l}$ has to be finite and simple. Hence it is of the form  $(C_p,\l_{C_p})$. But   $(M,\l_{M})$ is non-degenerate if  additionally we have $$1<n({\Z}\triangleright M)=\vert Hom_{Rings/{F_1}}({\Z}[x],{\Z}\triangleright M)\vert < \infty.$$
So we have to check when the number of $\L$-ring maps ${\Z}[x]\to {\Z}\triangleright M$ over $F_1$ is finite but bigger $1$ (we always have a zero morphism), where  $(M,\l_{M})$ is of the form $(C_p,\l_{C_p})$. Observe that if $\varphi \in Hom_{Rings/{F_1}}({\Z}[x],{\Z}\triangleright M)$ then $\varphi(x)=(0,m)$ for a certain $m\in M$.
Because $\varphi$ is a $\L$-ring homomorphism it has to commute with $\L$-operations on the source and the target. Recall that in ${\Z}[x]$, $\l^n(x)=0$ for $n>1$. Hence for $n>1$ we calculate
\bigskip

\noindent{\bf 4.1.1}$$0=\varphi (\l^n(x))=\l^n(\varphi(x))=\l^n((0,m)).$$  It means that we are looking for such objects  $(C_p,\l_{C_p})$ which give us vanishing of higher $\l$-operations on elements $(0,m)\in {\Z}\triangleright C_p$.  Observe that in general $\l_M$ is related to the $\L$-operations on ${\Z}\triangleright M$ via the Artin-Hasse isomorphism so we have to check that formulas $\l^n((0,m))=0$ imply $\l_{C_p,n}=0$ for $n>1$.

We are going to use now  the calculation from [H, Addendum 2.3], where the relation between the sequence $\l_{M,n}$ and $\l_{R\triangleright M}$ is calculated in full generality. We get that

$$\l_{{\Z}\triangleright C_p}((0,m))=((0,\l_{C_p,1}(m)),(0,\l_{C_p,2}(m)),(0,\l_{C_p,3}(m)),...)$$

Now we can use the general observation about the Artin-Hasse invariant. If $f(t)=1+\Sigma a_it^i\in \L(R)$ then we write $f(t)=\prod(1-(-1)^ib_it^i)$ and the Artin-Hasse isomorphism $E:\L(R)\to W(R)$ takes $f$ to the sequence $(b_1,b_2,b_3,...)$. Observe that if

$\ \ \ \ \ \ \ \ \ ({\bf *})\ \ \ \ b_i\cdot b_j=0\ \  {\rm for}\ \  {\rm any}\ \  i\ \  {\rm and}\ \  j$
\medskip
\noindent then up to sign $(b_1,b_2,b_3,...)=(a_1,a_2,a_3,...)$.

 Coming back to 4.4.1 we get that  up to sign:
 $$0=\l^n((0,m))=(0,\l_{C_p,n}(m))$$
  This implies that in our case $\l_{C_p,n}=0$ for $n>1$, as desired.
\bigskip

{\bf Corollary 4.5:} Recall that $\zeta_R$ denotes the Riemann $\zeta$-function of integers. We have
$$\zeta_R =\zeta_g(s,F_1-mod^{\l})$$

\bigskip
FINAL REMARK. It seems that the success of the Deligne's approach to Weil conjectures was caused by the fact that for an affine variety over a finite field we have two ways of calculations the $\zeta$ function. One is classical (categorical) via looking at finite simple objects in corresponding category. The second is via algebraic closure, rich geometric structure  and calculation of the same function via counting fixed points of the Frobenius action. In the present paper we tried to justify the statement that in the category of $\L$-rings the same two approaches should work.

\vfill
\eject

{\bf V. Bibliography:}
\bigskip

\noindent [Bec] - J. M. Beck: {\it Triples, algebras and cohomology}. Reprints in Theory and Applications of Categories 2 (2003) 1-59.

\noindent [Be] - S. Betley: {\it Riemann zeta via the category of monoids}. Kodai Math. Journal 36 (3) (2013) 487-490.

\noindent [B] - J. Borger: {\it $\L$-rings and the field with one element}. Preprint, 2009.

\noindent [BS] - J. Borger, B. de Smit: {\it Galois theory and integral models of $\L$-rings}. Bull. London Math. Soc. 40 (2008) 439-446.

\noindent [C] - A. Connes, C. Consani, M. Marcolli: {\it Fun with
$F_1$}. J. Number Theory 129 (2009) 1532-1561.

\noindent [D] - A. Deitmar: {\it Schemes over $F_1$}. Progr. Math
239 (2005) 87-100.

\noindent [H] - L. Hesselholt: {\it The big de Rham-Witt complex}. Preprint 2010.

\noindent [K] - N.Kurokawa: {\it Zeta Functions of Categories}.
Proc. Japan Acad. 72 (1996) 221-222.

\noindent [KOW] - N.Kurokawa, H. Ochiai, M.Wakayama: {\it Absolute
derivations and zeta functions}. Doc. Math., Extra Volume for
Kazuya Kato's 50 Birthday, (2003) 565-584.

\noindent [S] - C.Soule: {\it Les vari\'{e}t\'{e}s sur le corps a un
\'{e}l\'{e}ment}. Mosc. Math. J. 4 (2004), 217--244.

\noindent [Y] - D. Yau: {\it Lambda-Rings}. World Scientific 2010.
\bigskip

\bf Instytut Matematyki, University of Warsaw

ul.Banacha 2, 02-097 Warsaw, Poland

e-mail: betley@mimuw.edu.pl

\bye